\documentclass[reqno]{amsart}
\usepackage{hyperref}

\def\N{{\mathbb{N}}}

\def\R{{\mathbb{R}}}

\begin{document}
\title[Multiplier Rules]{On the multiplier rules}

\author[Jo\"el Blot]
{Jo${\rm \ddot e}$l Blot} 

\address{Jo\"{e}l Blot: Laboratoire SAMM EA4543,\newline
Universit\'{e} Paris 1 Panth\'{e}on-Sorbonne, centre P.M.F.,\newline
90 rue de Tolbiac, 75634 Paris cedex 13,
France.}
\email{blot@univ-paris1.fr}
\date{September 3, 2014}

\numberwithin{equation}{section}
\newtheorem{theorem}{Theorem}[section]
\newtheorem{lemma}[theorem]{Lemma}
\newtheorem{example}[theorem]{Example}
\newtheorem{remark}[theorem]{Remark}
\newtheorem{definition}[theorem]{Definition}
\newtheorem{proposition}[theorem]{Proposition}
\newtheorem{corollary}[theorem]{Corollary}
\begin{abstract}
We establish new results of first-order necessary conditions of optimality for finite-dimensional problems with inequality constraints and for problems with equality and inequality constraints, in the form of John's theorem and in the form of Karush-Kuhn-Tucker's theorem. In comparison with existing results we weaken assumptions of continuity and of differentiability.  
\end{abstract}
\maketitle
M.S.C. 2010: 90C30, 49K99.
\section{Introduction}
We consider first-order necessary conditions of optimality for finite-dimensional problems under inequality constraints and under equality and inequality constraints.
\vskip1mm
Let $\Omega$ be a nonempty open subset of $\R^n$, let $f_i : \Omega \rightarrow \R$ (when $i \in \{ 0,...,m \}$) be functions, let $\phi : \Omega \rightarrow \R$, $g_i : \Omega \rightarrow \R$ (when $i \in \{ 1,...,p \}$) and $h_j : \Omega \rightarrow \R$ (when $j \in \{1,...,q \}$) be functions. With these elements, we build the two following problems:
\[
({\mathcal I})
\left\{
\begin{array}{rl}
{\rm Maximize} & f_0(x)\\
{\rm when}& x \in \Omega\\
{\rm and} \;\; {\rm when} & \forall i \in \{1,...,m \}, f_i(x) \geq 0,
\end{array}
\right.
\]
and 
\[
({\mathcal M})
\left\{
\begin{array}{rl}
{\rm Maximize} & \phi(x)\\
{\rm when}& x \in \Omega\\
{\rm when} & \forall i \in \{1,...,p \}, g_i(x) \geq 0\\
{\rm and} \;\; {\rm when} & \forall j \in \{1,...,q \}, h_j(x) = 0.
\end{array}
\right.
\]
\vskip1mm
We provide necessary conditions of optimality under the form of Fritz John's conditions and under the form of Karush-Kuhn-Tucker's conditions. Our aim is to weaken the assumptions which permit to obtain such results. We can delete certain conditions of continuity and we can replace certain conditions of Fr\'echet-differentiability by conditions of G\^ateaux-differentiability.
\vskip1mm
The Farkas-Minkowski Theorem is one of the main tools that we use to establish the result for the  problem $({\mathcal I})$ with inequality constraints. A (local) theorem of implicit function permits to transform (locally) a problem with equality and inequality constraints (like $({\mathcal M})$) into a problem with only inequality constraints (like $({\mathcal I})$); it is why the Implicit Function Theorem of Halkin is one of the main tools to establish our result for $({\mathcal M})$.
\vskip3mm
These results are usual when we assume that all the functions are continuously Fr\'echet-differentiable on a neighborhood of $\hat{x}$, (\cite{ATF}, Chapter 3, Scetion 3.2)), (\cite{PB}, Chapitre 13, Section 2), \cite{HU}, \cite{BFS}, \cite{St}, \cite{MF}, (\cite{Tr}, Chapter 11). In \cite{Ha}, Halkin gives a multiplier rule only using the continuity on a neighborhood of $\hat{x}$ and the Fr\'echet-differentiability at $\hat{x}$ of the functions. His proof uses his implicit function theorem (Theorem \ref{th22}). In (\cite{Mi}, Chapter 24, Section 24.7) Michel gives another proof of the result of Halkin without to use an implicit function theorem but nevertheless using the Fixed Point Theorem of Brouwer. The proof of Michel is also explained in (\cite{BH}, Appendix B). In \cite{Po} we find a result for (${\mathcal I}$) with only the Fr\'echet-differentiability of the functions $f_i$ at $\hat{x}$.
\vskip1mm
There exist several works on the multiplier rules for locally Lipschitzian functions which are obtained by using the Clarke Calculus \cite{Cl}. After a famous theorem of Rademacher on the Lebesgue-almost everywhere Fr\'echet-differentiability of a locally Lipschtzian mapping, and since the Clarke-gradient is a upper semicontinuous correspondence, we can say that the locally Lipschitzian generalize the continuously Fr\'echet-differentiable mappings. Note that a mapping which is only Fr\'echet-differentiable (even all over a naighborhood of a point) is not necessarily locally Lipschtzian and a locally Lipschitzian mapping is not necessarily Fr\'echet-differentiable at a given point. And so there exist two different ways for the generalisation of the multiplier rules of the continuously differentiable setting: the locally Lipschitzian setting, and the (only) Fr\'echet-differentiable (or differentiability in a weaker sense than this one of Fr\'echet) setting. Our paper belongs to the second way.
\vskip1mm
Now we briefly describe the contents of the paper. In Section 2 we precise our notation and we recall two important tools. In Section 3 we state the new results for $({\mathcal I})$ and for $({\mathcal M})$. In Section 4 we prove the theorem of necessary condition of optimality for $({\mathcal I})$, and in Section 5, we prove the theorem of necessary condition of optimality for $({\mathcal M})$.
\section{Notation and recall}
First we precise the used notions of differentiability. Let $E$ and $F$ be two real normed spaces, let $\Omega$ be a nonempty open subset of $E$, $f : \Omega \rightarrow F$ be a mapping and let $x \in \Omega$ and $v \in E$. When it exists, the directional derivative of $f$ at $x$ in the direction of $v$ is $\vec{D} f(x; v) := \frac{d}{dt}_{\vert_{t=0}} f(x + tv)$. When $\vec{D} f(x; v)$ exists for all $v \in E$ and when $[v \mapsto \vec{D} f(x; v)]$ is linear continuous, we say that $f$ is G\^ateaux-differentiable at $x$; its G\^ateaux-differential at $x$ is $D_Gf(x) \in {\mathcal L}(E,F)$ (the vector space of the linear continuous mappings from $E$ into $F$) defined by $D_G f(x).v := \vec{D} f(x; v)$. The mapping $f$ is Fr\'echet-differentiable at $x$ when there exists $Df(x) \in {\mathcal L}(E,F)$ (so-called the Fr\'echet-differential of $f$ at $x$) and a mapping $\rho : \Omega - x \rightarrow F$ such that $\lim_{v \rightarrow 0} \rho (v) = 0$ and $f(x+v) = f(x) + Df(x).v + \Vert v \Vert \rho(v)$ for all $v \in \Omega -x$. When $f$ is Fr\'echet-differentiable at $x$ then $f$ is G\^ateaux-differentiable at $x$, and $D_Gf(x) = Df(x)$. When $E = E_1 \times E_2$, when $k \in \{ 1,2 \}$, $D_k f(x)$ (respectively $D_{G,k} f(x)$) denotes the partial Fr\'echet (respectively G\^ateaux)-differential of $f$ at $x$ with respect to the $k$-th variable. For all these notions we refer to the books (\cite{ATF}, Chapter 2, Section 2.2) and (\cite{Fl}, Chapter 4, sections 4.1, 4.2).
\vskip1mm
$\N$ denotes the set of the non negative integer numbers, $\N_* : \N \setminus \{ 0 \}$, $\R$ denotes the set of the real numbers and $\R_+$ denotes the set of the non negative real numbers. When $n \in \N_*$, we write $\R^{n*} := {\mathcal L}(\R^n, \R)$ the dual space.
\vskip1mm
We recall the Farkas-Minkowski Theorem.
\vskip1mm
\begin{theorem}\label{th21}
Let $m, n \in \N_*$, $\varphi_1$, ..., $\varphi_m \in \R^{n*}$, and $a \in \R^{n*}$. The two following assertions are equivalent.
\begin{enumerate}
\item[(i)] For all $x \in \R^n$, $(\forall i \in \{i,...,m \}, \varphi_i.x \geq 0) \Longrightarrow (a.x \geq 0)$.
\item[(ii)] There exists $\lambda^1$, ..., $\lambda^m \in \R_+$ such that $a = \sum_{1 \leq i \leq m} \lambda^i \varphi_i$.
\end{enumerate}
\end{theorem}
\vskip2mm
A complete proof of this result is given in (\cite{Ti}, Chapter, Sections 4.14-4.19) and in (\cite{Ka}, Chapter 2, Sections 2.5, 2.6). This result is present in many books like, for example (\cite{PB}, Chapter 13, Section 2), (\cite{Tr}, p. 176), (\cite{Be}, p. 164). A main difficulty of the proof of this theorem is the closedness of a finitely generated convex cone; a difficulty which is not ever well treated.
\vskip1mm
A second fundamental tool that we recall is the Implicit Function Theorem of Halkin for the Fr\'echet-differentiable mappings which are not necessarily continuously Fr\'echet-differentiable.
\vskip2mm
\begin{theorem}\label{th22}
Let $X$, $Y$, $Z$ be three real finite-dimensional mormed vector spaces, let $A \subset X \times Y$ be a nonempty open subset, let $f : A \rightarrow Z$ be a mapping, and let $(\bar{x}, \bar{y}) \in A$. We assume that the following conditions are fulfilled.
\begin{enumerate}
\item[(i)] $f(\bar{x}, \bar{y}) = 0$.
\item[(ii)] $f$ is continuous on a neighborhood of $(\bar{x}, \bar{y})$
\item[(iii)] $f$ is Fr\'echet-differentiable at $V$ and the partial Fr\'echet-differential $D_2f(\bar{x}, \bar{y})$ is bijective.
\end{enumerate}
Then there exist a neighborhood $U$ of $\bar{x}$ in $X$, a neighborhood $V$ of $\bar{y}$ in $Y$ such that $U \times V \subset A$, and a mapping $\psi : U \rightarrow V$ which satisfy the following conditions.
\begin{enumerate}
\item[(a)] $\psi(\bar{x}) = \bar{y}$
\item[(b)] For all $x \in U$, $f(x, \psi(x)) = 0$
\item[(c)] $\psi$ is Fr\'echet-differentiable at $\bar{x}$ and $D\psi(\bar{x}) = - D_2f(\bar{x}, \bar{y})^{-1} \circ D_1f(\bar{x}, \bar{y})$.
\end{enumerate}
\end{theorem}
\vskip2mm
This result is proven in \cite{Ha}. Its proof uses the Fixed Point Theorem of Brouwer. The electronic paper of Border \cite{Bo} is very useful to understand the role of each assumption of the theorem. Halkin does not use an open subset $A$; his function is defined on $X \times Y$. But it is easy to adapt his result. Since $\psi$ is Fr\'echet-differentiable at $\bar{x}$, $\psi$ is continuous at $\bar{x}$ and then we can consider a neighborhood of $\bar{y}$ and a neighborhood $U$ of $\bar{x}$ such that $\psi(U) \subset V$ and such that $U \times V \subset A$. 
\section{The main results}
For the problem (${\mathcal I}$) we state the following result.
\begin{theorem}\label{th31}
Let $\hat{x}$ be a solution of (${\mathcal I}$). We assume that the following assumptions are fulfilled.
\begin{enumerate}
\item[(i)] For all $i \in \{ 1,...,m \}$, $f_i$ is G\^ateaux-differentiable at $\hat{x}$.
\item[(ii)] For all $i \in \{ 1,...,m \}$, $f_i$ is lower semicontinuous at $\hat{x}$ when $f_i(\hat{x}) > 0$.
\end{enumerate}
Then there exist $\lambda^0$,..., $\lambda^m \in \R_+$ such that the following conditions hold.
\begin{enumerate}
\item[(a)] $(\lambda^0,..., \lambda^m) \neq (0,..., 0)$.
\item[(b)] For all $i \in \{ 1,..., m \}$, $\lambda^i f_i(\hat{x}) = 0$.
\item[(c)] $\sum_{0 \leq i \leq m} \lambda^i D_G f_i(\hat{x}) = 0$.
\end{enumerate}
\vskip1mm
If, in addition, we assume that the following assumption is fulfilled,
\begin{enumerate}
\item[(iii)] There exists $w \in \R^n$ such that, for all $i \in \{ 1,...,m \}$, $D_Gf_i(\hat{x}).w > 0$ when $f_i(\hat{x}) = 0$,
\end{enumerate}
then we can take $\lambda^0 = 1$.
\end{theorem}
\vskip3mm
The notion of lower semicontinuity is the classical one; see for instance \cite{Be} (p.74). For the problem (${\mathcal M}$), we state the following result.
\vskip2mm
\begin{theorem}\label{th32}
Let $\hat{x}$ be a solution of (${\mathcal M}$). We assume that the following assumptions are fulfilled.
\begin{enumerate}
\item[(i)] $\phi$ is Fr\'echet-differentiable at $\hat{x}$.
\item[(ii)] For all $i \in \{ 1,..., p \}$, $g_i$ is Fr\'echet-differentiable at $\hat{x}$ when $g_i(\hat{x}) = 0$.
\item[(iii)] For all $i \in \{ 1,..., p \}$, $g_i$ is lower semicontinuous at $\hat{x}$ and G\^ateaux-dif\-
ferentiable at $\hat{x}$ when $g_i(\hat{x}) > 0$.
\item[(iv)] For all $j \in \{ 1,..., q \}$, $h_j$ is continuous on a neighborhood of $\hat{x}$ and Fr\'echet-differentiable at $\hat{x}$.
\end{enumerate}
Then there exist $\lambda^0$, $\lambda^1$,..., $\lambda^p \in \R_+$ and $\mu^1$,..., $\mu^q \in \R$ such the following conditions are satisfied.
\begin{enumerate}
\item[(a)] $(\lambda^0, \lambda^1,..., \lambda^p, \mu^1, ..., \mu^q) \neq (, 0, ..., 0)$.
\item[(b)] for all $i \in \{ 1,..., p \}$, $\lambda^i g_i(\hat{x}) = 0$.
\item[(c)] $\lambda^0 D \phi(\hat{x}) + \sum_{1 \leq i \leq p} \lambda^i D_G g_i(\hat{x}) + \sum_{1 \leq j \leq q} \mu^j Dh_j(\hat{x}) = 0$.
\end{enumerate}
\vskip1mm
Moreover, under the additional assumption
\begin{enumerate}
\item[(v)] $Dh_1(\hat{x})$,..., $Dh_q(\hat{x})$ are linearly independent,
\end{enumerate}
we can take 
\begin{enumerate}
\item[(d)] $(\lambda^0, \lambda^1,..., \lambda^p) \neq (0, 0, ..., 0)$.
\end{enumerate}
\vskip1mm
Furthermore, under (v) and under the additional assumption
\begin{enumerate}
\item[(vi)] There exists $w \in \bigcap_{1 \leq j \leq q} {\rm Ker} D h_j(\hat{x})$ such that, for all $i \in \{ 1, ..., p \}$,\\
 $Dg_i(\hat{x}).w > 0$ when $g_i(\hat{x}) = 0$,
\end{enumerate}
we can take
\begin{enumerate}
\item[(e)] \hskip3mm $\lambda^0 = 1$.
\end{enumerate}
\end{theorem}
\begin{remark}
The assumption (iii) is generally called the Mangarasian-Fromowitz's condition. In (\cite{PB}, p. 289) the author associates this condition at a work of Abadie in 1965 (it is difficult to find the reference). In (\cite{Tr}, p. 197) we find a catalog of the variations of this condition due to Cottle, Zandwill, Kuhn and Tucker, and Abadie.
\end{remark}
\vskip1mm
In comparison with the Halkin's multiplier rule, for problem (${\mathcal I}$) we have deleted the assumptions of local continuity on a neighborhood of $\hat{x}$ of the $f_i$ and we have replaced their Fr\'echet-differentiability by their G\^ateaux-differentiability, and for problem (${\mathcal M}$), we have deleted the assumptions of local continuity on $\phi$ and on the $g_i$. In comparison with the result of \cite{Po} for problem (${\mathcal I}$), we have replaced the Fr\'echet-differentiability of the $f_i$ by their G\^ateaux-differentiability. Note that the G\^ateaux-differentiability of a mapping at a point does not imply the continuity of this mapping at this point.
\section{Proof of Theorem \ref{th31}}
Doing a change of index, we can assume that $\{ 1,...,e \} = \{ i \in \{ 1,...,m \} : f_i(\hat{x}) = 0 \}$. If $f_i(\hat{x}) > 0$ for all $i \in \{ 1,...,m \}$, the using the lower semicontinuity of (ii), there exists an open neighborhood of $\hat{x}$ on which $\hat{x}$ maximizes $f_0$ (without constraints). Then using (i) we obtain $D_Gf_0(\hat{x}) = 0$, and we conclude by taking $\lambda^0 := 1$ and $\lambda^i := 0$ for all $i \in \{ 1,...,m \}$. And so, for the sequel of the proof we assume that $1 \leq e \leq p$.
\subsection{Proof of (a), (b), (c)}
Ever using (ii) we can assert that there exists an open neighborhood $\Omega_1 \subset \Omega$ of $\hat{x}$ such that, for all $x \in \Omega_1$ and for all $i \in \{e+1,...,m \}$, $f_i(x) > 0$ when $e < m$. When $e = m$ we simply take $\Omega_1 := \Omega$. Then $\hat{x}$ is a solution of the following problem.
\[
({\mathcal P})
\left\{
\begin{array}{rl}
{\rm Maximize} & f_0(x)\\
{\rm when}& x \in \Omega_1\\
{\rm and} \; \; {\rm when} & \forall i \in \{1,...,e \}, \; f_i(x) \geq 0.
\end{array}
\right.
\]
For all $k \in \{ 0,...,e \}$ we introduce the set
\begin{equation}\label{eq41}
A_k := \{ v \in \R^n : \forall i \in \{k,...,e \}, D_Gf_i(\hat{x}).v > 0 \}.
\end{equation}
We want to prove that $A_0 = \emptyset$. To realize that, we proceed by contradiction; we assume that $A_0 \neq \emptyset$, and so there exists $w \in \R^n$ such that 
$D_Gf_i(\hat{x}).w > 0$ for all $i \in \{0,...,e \}$. Since $\Omega_1$ is open, there exists $\theta_* \in (0, + \infty)$ such that $\hat{x} + \theta w \in \Omega_1$ for all $\theta \in [0, \theta_*]$. After (i), for all $i \in \{0,..., e \}$, the function $\sigma_i : [0, \theta_*] \rightarrow \R$, defined by $\sigma_i( \theta) := f_i(\hat{x} + \theta w)$, is differentiable at $0$, and its derivative is $\sigma_i'(0) = D_Gf_i(\hat{x}).w$. The differentiability of $\sigma_i$ at $0$ implies the existence of a function $\rho_i : [0, \theta_*] \rightarrow \R$ such that $\lim_{\theta \rightarrow 0} \rho_i(\theta) = 0$ and such that $\sigma_i(\theta) = \sigma_i(0) + \sigma_i'(0) \theta + \theta \rho_i(\theta)$ for all $\theta \in [0, \theta_*]$. Translating this last equality we obtain $f_i(\hat{x} + \theta w) = f_i(\hat{x}) + \theta (D_Gf_i(\hat{x}).w + \rho_i(\theta))$. Since $D_Gf_i(\hat{x}).w > 0$ and since $\lim_{\theta \rightarrow 0} \rho_i(\theta) = 0$, we obtain the existence of $\theta_i \in (0, _theta_*]$ such that $D_Gf_i(\hat{x}).w + \rho_i(\theta) > 0$ for all $\theta \in (0, \theta_i]$. Setting $\hat{\theta} := \min \{ \theta_i : i \in \{0,...,e \}$ we obtain that $f_i(\hat{x} + \theta w) > f_i((\hat{x})$ for all $\theta \in (0, \hat{\theta}]$ and for all $i \in \{0,...,e \}$. Then using $i \in \{ 1,...,e \}$, this last relation ensures that $\hat{x} + \theta w$ is admissible for (${\mathcal P}$) when $\theta \in (0, \hat{\theta}]$, and using this last relation when $i=0$ we obtain $f_0(\hat{x} + \theta w) > f_0((\hat{x})$ when $\theta \in (0, \hat{\theta}]$, that is impossible since $\hat{x}$ is a solution of (${\mathcal P}$). And so the reasoning by contradiction is complete, and we have proven
\begin{equation}\label{eq42}
A_0 = \emptyset.
\end{equation}
When $A_e = \emptyset$ there is not any $v \in \R^n$ such that $D_Gf_e(\hat{x}).v > 0$, that implies that $D_Gf_e(\hat{x}) = 0$. Then taking $\lambda^e := 1$ and $\lambda^i := 0$ when $i \in \{0,...,m \} \setminus \{e \}$, we obtain the conclusions (a), (b), (c). And so we have proven
\begin{equation}\label{eq43}
A_e = \emptyset \Longrightarrow ((a), (b), (c) \;\; {\rm hold}).
\end{equation}
Now we assume that $A_e \neq \emptyset$. Since we have $A_0 = \emptyset$ after (\ref{eq42}) and $A_i \subset A_{i+1}$ we can define
\begin{equation}\label{eq44}
k := \min \{ i \in \{1,..., e \} : A_i \neq \emptyset \}.
\end{equation}
Note that $A_k \neq \emptyset$ and that $A_{k-1} = \emptyset$. We consider the following problem
\[
({\mathcal Q}) 
\left\{
\begin{array}{rl}
{\rm Maximize} & D_Gf_{k-1}(\hat{x}).v\\
{\rm when} & v \in \R^n\\
{\rm and} \;\; {\rm when} & \forall i \in \{ k,...,e \}, \;\; D_Gf_i(\hat{x}).v \geq 0.
\end{array}
\right.
\]
We want to prove that $0$ is a solution of (${\mathcal Q}$). To do that, we proceed by contradiction; we assume that there exists $y \in \R^n$ such that ($\forall i \in \{k,...,e \}$, $D_Gf_i(\hat{x}).y \geq 0$) and $D_Gf_{k-1}(\hat{x}).y > 0 = D_Gf_{k-1}(\hat{x}).0$. Since $A_k \neq \emptyset$, there exists $z \in \R^n$ such that $D_Gf_i(\hat{x}).z > 0$ when $i \in \{k,...,e \}$. We cannot have $D_Gf_{k-1}(\hat{x}).z > 0$ since $A_{k-1} = \emptyset$. Therefore we have $D_Gf_{k-1}(\hat{x}).z \leq 0$. If $D_Gf_{k-1}(\hat{x}).z < 0$ we choose $\epsilon$ such that $0< \epsilon < \frac{D_Gf_{k-1}(\hat{x}).y}{D_Gf_{k-1}(\hat{x}).z}$. Then we have $D_Gf_{k-1}(\hat{x}).y + \epsilon D_Gf_{k-1}(\hat{x}).z > 0$. If $D_Gf_{k-1}(\hat{x}).z = 0$ we arbitrarily choose $\epsilon \in (0, + \infty)$ and we have also $D_Gf_{k-1}(\hat{x}).y + \epsilon D_Gf_{k-1}(\hat{x}).z > 0$. We set $u_{\epsilon} := y + \epsilon z$, and we note that $D_Gf_{k-1}(\hat{x}).u_{\epsilon} = D_Gf_{k-1}(\hat{x}).y + \epsilon D_Gf_{k-1}(\hat{x}).z > 0$. Furthermore, when $i \in \{ k,...,e \}$, we have $D_Gf_i(\hat{x}).u_{\epsilon} = D_Gf_{i}(\hat{x}).y + \epsilon D_Gf_{i}(\hat{x}).z > 0$ since the three terms are positive. Therefore we have $u_{\epsilon} \in A_{k-1}$ that is impossible since $A_{k-1} = \emptyset$. And so the reasoning by contradiction is complete, and we have proven
\begin{equation}\label{eq45}
A_e \neq \emptyset \Longrightarrow (0 \;\; {\rm solves} \;\; ({\mathcal Q})).
\end{equation}
Since $0$ solves (${\mathcal Q}$), we have, for all $v \in \R^n$, 
$$(\forall i \in \{ k,...,e \}, \; \; D_Gf_i(\hat{x}).v \geq 0) \Longrightarrow (D_Gf_{k-1}(\hat{x}).v \geq 0).$$
Then we use Theorem \ref{th21} that ensures the existence of $\alpha^k$,..., $\alpha^e \in \R_+$ such that $D_Gf_{k-1}(\hat{x}) + \sum_{k \leq i \leq e} \alpha^i D_Gf_{i}(\hat{x}) = 0$. We set
\[
\lambda^i := \left\{
\begin{array}{ccl}
0 & {\rm if} & i \in \{0,...,k-2 \}\\
1 & {\rm if} & i= k-1\\
\alpha^i & {\rm if} & i \in \{k,...,e \}\\
0 & {\rm if} & i \in \{e+1, ..., m \},
\end{array}
\right.
\]
and we obtain
\begin{equation}\label{eq46}
A_e \neq \emptyset \Longrightarrow ((a), (b), (c) \;\; {\rm hold}).
\end{equation}
Then, with (\ref{eq43}) and (\ref{eq46}) the conclusions (a), (b), (c) are proven.
\subsection{Proof of (d)}
The assumption (iii) means that $A_1 \neq \emptyset$, and by (\ref{eq42}) we know that $A_0 = \emptyset$. Proceeding like in the proof of (\ref{eq45}) we prove that $0$ is a solution of the following problem
\[
\left\{
\begin{array}{rl}
{\rm Maximize} & D_Gf_0(\hat{x}).v\\
{\rm when} & v \in \R^n\\
{\rm and} \;\; {\rm when} & \forall i \in \{1,...,e \}, \;\; D_Gf_i(\hat{x}).v \geq 0.
\end{array}
\right.
\]
Then using Theorem \ref{th21}, there exist $\alpha^1$,..., $\alphaê \in \R_+$ such that 
$$D_Gf_0(\hat{x}) + \sum_{1 \leq i \leq e} \alpha^i D_Gf_i(\hat{x})= 0.$$
We conclude by setting
\[
\lambda^i := 
\left\{
\begin{array}{ccl}
1 & {\rm if} & i=0\\
\alpha^i & {\rm if} & i \in \{ 1,..., e \}\\
0 & {\rm if} & i \in \{ e+1,..., m \}.
\end{array}
\right.
\]
And so the proof of Theorem \ref{th31} is complete.
\begin{remark}
The use of the sets $A_k$ comes from the book of Alexeev-Tihomirov-Fomin \cite{ATF}, end the proof of formula (\ref{eq46}) is similar to their proof (p. 247-248). The use of the set $A_0$ is yet done in \cite{Ha}.
\end{remark}

\section{Proof of Theorem \ref{th32}}
We split this proof in seven steps.
\subsection{First step : a first simple case.}
If $Dh_1(\hat{x})$, ..., $Dh_q(\hat{x})$ are linearly dependent, there exist $\mu^1$,..., $\mu^q \in \R$ such that $(\mu^1, ..., \mu^q) \neq (0,...,0)$ and such that $\sum_{1 \leq j \leq q} \mu^j Dh_j(\hat{x}) = 0$. Then it suffices to take $\lambda^i = 0$ for all $i \in \{0,...,p \}$ to obtain the conclusions (a), (b), (c).
\vskip1mm
Now in the sequel of the proof {\bf we assume that the assumption (v) is fulfilled}.
\subsection{Second step : To delete the non satured inequality constraints.}
Doing a change of index, we can assume that $\{1,...,e \} := \{ i \in \{1,...,p \} : g_i(\hat{x}) = 0 \}$. Using the lower semicontinuity at $\hat{x}$ of the $g_i$ when $i \in \{ e+1,...,p \}$, we can say that there exists an open neighborhood $\Omega_1$ of $\hat{x}$ in $\Omega$ such that $g_i(x) > 0$ when $x \in \Omega_1$ and when $i \in \{ e+1,...,p \}$. And so $\hat{x}$ is a solution of the following problem
\[
({\mathcal M}_1)
\left\{
\begin{array}{rl}
{\rm Maximize} & \phi(x)\\
{\rm when} & x \in \Omega_1\\
{\rm when} & \forall i \in \{1,...,e \}, \;\; g_i(x) \geq 0\\
{\rm and} \;\; {\rm when} & \forall j \in \{1,...,q \}, \;\; h_j(x) = 0.
\end{array}
\right.
\]
\subsection{To delete the equality constraints.}
We consider the mapping $h : \Omega_1 \rightarrow \R^q$ defined by $h(x) := (h_1(x),...,h_q(x))$. Under (iv) and (v), $h$ continuous on a neighborhood of $\hat{x}$, and it is Fr\'echet-differentiable at $\hat{x}$ with $Dh(\hat{x})$ onto.
\vskip1mm
We set $E_1 := {\rm Ker} Dh(\hat{x})$ and we take a vector subspace of $\R^n$ such that $E_1 \oplus E_2 = \R^n$. And we can do the assimilitation $\R^n = E_1 \times E_2$. We set $(\hat{x}_1, \hat{x}_2) := \hat{x} \in E_1 \times E_2$. Then the partial differential $D_2h(\hat{x})$ is an isomorphism from $E_2$ onto $\R^q$. Now we can use Theorem \ref{th22} and assert that there exist a neighborhood $U_1$ of $\hat{x}_1$ in $E_1$, a neighborhood $U_2$ of $\hat{x}_2$ in $E_2$, and a mapping $\psi : U_1 \rightarrow U_2$ such that $\psi(\hat{x}_1) = \hat{x}_2$, $h(x_1, \psi(x_1)) = 0$ for all $x_1 \in U_1$, and such that $\psi$ is Fr\'echet-differentiable at $\hat{x}_1$ with $D\psi(\hat{x}_1) = - D_2h(\hat{x})^{-1} \circ D_1h(\hat{x}) = 0$ since $D_1h(\hat{x}) = Dh(\hat{x})_{\vert_{E_1}}= 0$.
\vskip1mm
We define $f_0 : U_1 \rightarrow \R$ by setting $f_0(x_1) := \phi(x_1, \psi(x_1))$, and $f_i : U_1 \rightarrow \R$ by setting $f_i(x_1) := g_i(x_1, \psi(x_1))$ for all $i \in \{1,...,e \}$. Since $\hat{x}$ is a solution of (${\mathcal M}_1$), $\hat{x}_1$ is a solution of the following problem without equality constraints
\[
({\mathcal R})
\left\{
\begin{array}{rl}
{\rm Maximize} & f_0(x_1)\\
{\rm when} & x_1 \in U_1\\
{\rm and} \;\; {\rm when} & \forall i \in \{1,...,e \}; \;\; f_i(x_1) \geq 0.
\end{array}
\right.
\]
\subsection{Fourth step : To use Theorem \ref{th31}.}
Since $\psi$ is Fr\'echet-differentiable at $\hat{x}_1$, the mapping $[x_1 \mapsto (x_1, \psi(x_1))]$ is Fr\'echet-differentiable at $\hat{x}_1$, and using (i) and (ii), we obtain that $f_i$ is Fr\'echet-differentiable (and therefore G\^ateaux-differentiable) at $\hat{x}_1$, for all $i \in \{ 0,...,e \}$. Note that $f_i(\hat{x}_1) = 0$ for all $i \in \{ 1,...,e \}$. Consequently we can use Theorem \ref{th31} on (${\mathcal R}$) that permits us to ensure the existence of $\lambda^0$,$\lambda^1$,..., $\lambda^e \in \R_+$ such that
\begin{equation}\label{eq51}
(\lambda^0, \lambda^1,..., \lambda^e) \neq (0,0,...,0)
\end{equation}
\begin{equation}\label{eq52}
\forall i \in \{1,...,e \}, \;\; \lambda^i f_i(\hat{x}_1) = 0
\end{equation}
\begin{equation}\label{eq53}
\sum_{0 \leq i \leq e} \lambda^i D_Gf_i(\hat{x}_1) = 0.
\end{equation}
\subsection{The proof of (a), (b), (c).}
Since $D_Gf_0(\hat{x}_1) = Df_0(\hat{x}_1) = D_1 \phi(\hat{x}) + D_2 \phi(\hat{x}) \circ D \psi((\hat{x}_1) = D_1 \phi(\hat{x})$ since $D \psi(\hat{x}_1) = 0$,$ D_Gf_i(\hat{x}_1) = Df_i(\hat{x}_1) = D_1g_i(\hat{x}) + D_2 g_i(\hat{x}) \circ D \psi(\hat{x}_1) = D_1 g_i(\hat{x})$, for all $i \in \{1,...,e \}$, the formula (\ref{eq53}) implies
\begin{equation}\label{eq54}
\lambda^0 D_1 \phi(\hat{x}) + \sum_{1 \leq i \leq e} \lambda^i D_1g_i(\hat{x}) = 0.
\end{equation}
We set
\begin{equation}\label{eq55}
M := - (\lambda^0 D_2 \phi(\hat{x}) + \sum_{1 \leq i \leq e} \lambda^i D_2g_i(\hat{x})) \circ D_2h(\hat{x})^{-1}) \in \R^{q*}.
\end{equation}
Then we have 
$$\lambda^0 D_2 \phi(\hat{x}) + \sum_{1 \leq i \leq e} \lambda^i D_2g_i(\hat{x}) + M \circ D_2h(\hat{x}) = 0.$$
Denoting by $\mu^1$,..., $\mu^q \in \R$ the coordinates of $M$ in the canonical basis of $\R^{q*}$, we obtain
\begin{equation}\label{eq56}
\lambda^0 D_2 \phi(\hat{x}) + \sum_{1 \leq i \leq e} \lambda^i D_2g_i(\hat{x}) + \sum_{1 \leq j \leq q} \mu^j D_2h_j(\hat{x}) = 0.
\end{equation}
Since $E_1 = {\rm Ker}Dh(\hat{x}) = \bigcap_{1 \leq j \leq q} {\rm Ker}Dh_j(\hat{x})$, we have $D_1h_j(\hat{x}) = Dh(\hat{x})_{\vert_{E_1}} = 0$ for all $j$, from (\ref{eq54}) we obtain
\begin{equation}\label{eq57}
\lambda^0 D_1 \phi(\hat{x}) + \sum_{1 \leq i \leq e} \lambda^i D_1g_i(\hat{x}) + \sum_{1 \leq j \leq q} \mu^j D_1 h_j(\hat{x})= 0.
\end{equation}
From (\ref{eq56}) and (\ref{eq57}) we obtain
\begin{equation}\label{eq58}
\lambda^0 D \phi(\hat{x}) + \sum_{1 \leq i \leq e} \lambda^i Dg_i(\hat{x}) + \sum_{1 \leq j \leq q} \mu^j D h_j(\hat{x})= 0.
\end{equation}
We set $\lambda^i := 0$ when $i \in \{ e+1, ..., p \}$, and so (\ref{eq58}) implies (c). With (\ref{eq51}) we obtain (a), and with (\ref{eq52}) we obtain (b). And so the proof of (a), (b), (c) is complete.
\subsection{The proof of (d).} The relation (\ref{eq53}) provides the conclusion (d).
\subsection{The proof of (e).} 
When $i \in \{ 1,...,e \}$, we have yet seen that $Df_i(\hat{x}_1) = D_1g_i(\hat{x} = Dg_i(\hat{x})_{\vert_{E_1}}$. And so the translation of the assumption (vi) gives
$$\exists w \in E_1 \;\; {\rm s.t.} \;\; \forall i \in \{ 1,...,e \}, \;\; Df_i(\hat{x}_1).w > 0.$$
That permits us to use the last assertion of Theorem \ref{th31} on (${\mathcal R}$) to ensure that we can choose $\lambda^0 = 1$.\\
Then the proof of Theorem \ref{th32} is complete.
\begin{remark}
We see in this proof that the assumption of Fr\'echet-diffferentiability of the $h_j$ is used to can apply the Implicit Function of Halkin. The assumption of Fr\'echet-diffferentiability of $\phi$ and of the $g_i$ for which the associated constraint is satured is used to obtain the differentiability when we compose them with $ h_j$ (to obtain the differentiability of the $f_i$). The Hadamard-differentiability is sufficient to do that, but in finite-dimensional spaces, the Hadamard-differentiability and the Fr\'echet-differentiability coincide (\cite{Fl}, p. 266).
\end{remark}

\end{document}